\documentclass{commat}

\usepackage{tikz-cd}

\DeclareMathOperator\AQ{KQ}
\DeclareMathOperator\C{C}
\DeclareMathOperator\E{E}
\DeclareMathOperator\EQ{EQ}
\DeclareMathOperator\ET{ET}
\DeclareMathOperator\GL{GL}
\DeclareMathOperator\GQ{GQ}
\DeclareMathOperator\G{G}
\DeclareMathOperator\I{I}
\DeclareMathOperator\K{K}
\DeclareMathOperator\KH{KH}
\DeclareMathOperator\KQ{KQ}
\DeclareMathOperator\Max{Max}
\DeclareMathOperator\SL{SL}
\DeclareMathOperator\SQ{SQ}

\DeclareMathOperator\T{T}
\DeclareMathOperator\w{W}
\newcommand\kk{\textrm{K}_1}
\newcommand\nk{\textrm{NK}_1}
\newcommand\sk{\textrm{SK}_1}
\newcommand{\lra}{\longrightarrow}
\newcommand{\mbb}{\mathbb}
\newcommand{\ol}{\overline} 
\newcommand{\ra}{\rightarrow}

\newtheorem{nt}[definition]{Notation}

\title{%
    Results on ${\kk}$ of general quadratic groups
    }

\author{%
    Rabeya Basu and Kuntal Chakraborty
    }

\affiliation{
    \address{Rabeya Basu --
    Indian Institute of
    Science Education and Research (IISER) Pune,  India
        }
    \email{%
    rabeya.basu@gmail.com, rbasu@iiserpune.ac.in
    }
    \address{Kuntal Chakraborty --
    Indian Institute of
    Science Education and Research (IISER) Pune,  India
        }
    \email{%
    kuntal.math@gmail.com
    }
    }

\abstract{%
     In the first part of this article, we discuss the relative cases of Quillen--Suslin's local-global principle for the general quadratic (Bak's unitary) groups, and its applications for the (relative) stable and unstable ${\kk}$-groups. 
    The second part is dedicated to the graded version of the local-global principle for the general quadratic groups and its application to deduce a result for Bass' nil groups.
    }



\VOLUME{32}
\YEAR{2024}
\NUMBER{1}
\firstpage{175}
\DOI{https://doi.org/10.46298/cm.9855}

\begin{paper}

\section{Introduction} 

In \cite{RB} and \cite{basu2}, the first author has discussed many results in classical K-theory for the absolute cases, related to the Serre's problem on projective modules. In this article, following the tricks used in \cite{BRK}, we are going to consider some problems for the relative cases; {\it viz} Quillen--Suslin's local-global principle of the transvection subgroups, ${\kk}$-stabilization and 
the structure of unstable ${\kk}$-groups of the general quadratic (Bak's unitary) groups over associative rings which are finite over the center. For previous results on these problems we refer to the works of Bak--Basu--Rao--Khanna for the local-global principle (L-G principle) in 
\cite{brk}, \cite{BBR}, \cite{BRK}, Bak--Petrov--Tang for ${\kk}$-stabilization in 
\cite{BTP}, and Bak--Harzat--Vavilov for solvability of unstable ${\kk}$-groups in \cite{Bak}, \cite{HV}, \cite{HR}.

For the linear case, the graded version of L-G principle was studied by Chouinard in \cite{C}, 
and by Gubeladze in \cite{gubel}, \cite{gubel2}. In \cite{BS}, the first author and M.K. Singh deduced an analog for the traditional classical groups. In the second part of the article, we have deduced an analog for the transvection subgroups of the 
general quadratic groups over graded rings. As an application, by using a recent result on Higman linearization due to V. Kepoeiko (cf. \!\cite{V}, \cite{V1}),
we could revisit a problem on absence of torsion in Bass' nil group of the general quadratic groups, and its graded analog. We refer to \cite{STI}, \cite{BL1}, \cite{WEL}, \cite{WEL1}, and \cite{RB} for previous results in this direction.

\section{Preliminaries} 
Let us recall some necessary definitions and the key lemmas.

\begin{definition} [cf. \!\cite{Bak1}]
	\label{form ring}	
	Let $R$ be an (not necessarily commutative) associative ring with identity, and
		with involution $-:R\rightarrow R$, $a\mapsto \ol{a}$. Let $\lambda\in C(R)$ = center of $R$ be an element
		with the property $\lambda\ol{\lambda}=1$. We define additive subgroups of $R$:
	$$\Lambda_{max}=\{a\in R\mid a=-\lambda \ol{a}\}~ \text{and }~ \Lambda_{min}=\{a-\lambda\ol{a}\mid a\in R\}.$$
	One checks that $\Lambda_{max}$ and  $\Lambda_{min}$ are closed under the conjugation operation
		$a\mapsto \ol{x}ax$ for any $x\in R$. A $\lambda$-{\it form parameter} on $R$ is an additive subgroup  $\Lambda$ of $R$ such that $\Lambda_{min}\subseteq \Lambda \subseteq \Lambda_{max}$, and 
		$\ol{x}\Lambda x\subseteq \Lambda$ for all $x\in R$. A pair $(R,\Lambda)$ is called a {\it form ring}.
\end{definition}

\begin{remark}
		For a from ring $(R,\Lambda)$, we can extend the involution $-:R\to R$ to  an involution $-:R[X]\to R[X]$ by setting $\ol{X}=X$.
\end{remark}

\begin{lemma}
	Let $(R,\Lambda,\lambda)$ be a form ring. Then $(R[X],\Lambda[X],\lambda)$ is also a form ring
	with respect to the involution obtained by extending, as in the previous remark.
\end{lemma}

\begin{proof} Let us consider the subgroups $\Lambda_{min}(R[X])$ and $\Lambda_{max}(R[X])$. It can be checked that $\Lambda_{min}(R[X])=\Lambda_{min}(R)[X]$, and $\Lambda_{max}(R[X])=\Lambda_{max}(R)[X]$. Hence it follows
$$\Lambda_{min}(R[X])\subseteq \Lambda[X]\subseteq \Lambda_{max}(R[X]).$$ To prove that it is closed under conjugation, 
let us consider elements $a(X)\in R[X]$ and $b(X)\in \Lambda[X]$. Using double induction on the degrees of $a(X)$ and $b(X)$, we will show that 
\[\ol{a(X)}b(X)a(X)\in \Lambda[X].\]
We first show that the statement is true for $\deg(b(X))=0$. Assume $b(X)=b\in \Lambda$. Hence we need to prove 
$\ol{a(X)}b a(X)\in \Lambda[X]$ for all $a(X)\in R[X]$. Clearly, this is obvious when $\deg(a(X))=0$.  For a linear polynomial $a(X)=a_0+a_1X$, 
\begin{align*}
	\ol{a(X)}b a(X) &= (\ol{a}_0+\ol{a}_1X)b(a_0+a_1X)\\
	&= \ol{a}_0ba_0+(\ol{a}_0b a_1+\ol{a}_1b a_0)X+\ol{a}_1b a_1 X^2 \\
	&= \ol{a}_0b a_0+\ol{(a_0+a_1)}b(a_0+a_1)X-(\ol{a}_0ba_0+\ol{a}_1ba_1)X+\ol{a}_1ba_1X^2.
\end{align*}
Hence the statement is true for $\deg(b(X))=0$ and $\deg(a(X))=1$. Assume the statement for $\deg(b(X))=0$ is true for all $a(X)$ with $\deg(a(X))<m$.  Suppose $\deg(a(X))=m$. Then, $a(X)=a_{m-1}(X)+a_mX^m$, where $a_{m-1}(X)$ is a polynomial of degree at most $m-1$, and $a_m\in R$. Hence, we have 

\begin{align*}
	\ol{a(X)}ba(X)=&(\ol{a_{m-1}(X)}+\ol{a}_mX^m)b(a_{m-1}(X)+a_mX^m)\\
	=&\ol{a_{m-1}(X)}ba_{m-1}(X)+\ol{a_{m-1}(X)}ba_mX^m+\ol{a}_mba_{m-1}(X)X^m+\ol{a}_mba_mX^{2m}.
\end{align*} 

Now \begin{gather*}
	\ol{a_{m-1}(X)}ba_m+ \ol{a}_m ba_{m-1}(X) =\\
	=\ol{(a_{m-1}(X)+a_m)}b(a_{m-1}(X)+a_m)-(\ol{a_{m-1}(X)}ba_{m-1}(X)+\ol{a}_mb a_m).
\end{gather*}
Hence by induction it follows that $\ol{a(X)}b a(X)\in \Lambda[X]$.
Assume the statement is true for $\deg(b(X))<n$. Let $b^{\prime}(X)$ be a polynomial in $\Lambda[X]$ of degree $n$. Then we can write $b^{\prime}(X)=b^{\prime \prime}(X)+b_nX^n$, where $b^{\prime\prime}(X)$ is a polynomial in $\Lambda[X]$ of degree at most $n-1$ and $b_n\in \Lambda$.  Hence by induction

\begin{align*}
	\ol{a(X)}b^{\prime}(X)a(X)=&\ol{a(X)}(b^{\prime \prime}(X)+b_nX^n)a(X)\\
	=&\ol{a(X)}b^{\prime \prime}(X)(a(X))+(\ol{a(X)}b_n a(X))X^n \in \Lambda[X].
 \qedhere
\end{align*}
\end{proof}

To define Bak's unitary group or the general quadratic group, we fix a (non-zero) central element $\lambda \in R$ with $\lambda\ol{\lambda}=1$, and 
then consider the form $$\psi_n= \begin{pmatrix} 0 &  {\I}_n \\ \lambda {{\I}}_n &0\end{pmatrix}.$$

\noindent{\bf Bak's Unitary or General Quadratic Groups ${\GQ}$:}  

$${\GQ}(2n, R,\lambda) ~ = ~ \{\sigma\in {\GL}(2n, R)\,|\, \ol{\sigma}\psi_n
\sigma=\psi_n\}.$$

\subsection*{Elementary Quadratic Matrices:} Let $\rho: \{1,2,\ldots,n\} \to \{n+1, n+2, \ldots, 2n \}$ be defined by $\rho(i)=n+i$. 
Let $e_{ij}$ be the matrix with $1$ in the $ij$-th position
and $0$'s elsewhere. For $a\in R$, and $1\leq i,j\leq n$, we define

\begin{center}
	$q\varepsilon_{ij}(a)={\I}_{2n}+ae_{ij}-\ol{a}e_{\rho(j)\rho(i)}$ for $i\neq j$,
	
	\[ qr_{ij}(a) = \left\{\begin{array}{ll}
		{\I}_{2n}+ ae_{i\rho(j)}-\lambda\ol{a}e_{j\rho(i)} & \text{for}~ i\neq j \\{\I}_{2n}+ae_{\rho(i)j} & \text{for}~ i=j,
	\end{array} \right. \]  \\
	\[ ql_{ij}(a) = \left\{\begin{array}{ll}
		{\I}_{2n}+ ae_{\rho(i)j}-\ol{\lambda}\ol{a}e_{\rho(j)i} & \text{for}~ i\neq j 
		\\{\I}_{2n}+ae_{\rho(i)j} & \text{for}~ i=j.
	\end{array} \right. \]  \\
\end{center}
(Note that for the second and third type of elementary matrices, if $i=j$, then
we get $a=-\lambda \ol{a}$, and hence it forces that $a\in\Lambda_{max}(R)$. One checks that these
above matrices belong to ${\GQ}(2n,R,\Lambda)$; cf. \!\cite{Bak1}.) 

\textbf{$n$-th Elementary Quadratic Group} ${\EQ}(2n,R,\Lambda)$: The subgroup generated
by $q\varepsilon_{ij}(a),qr_{ij}(a) \text{and } ql_{ij}(a)$, for $a\in R$ and $1\leq i,j\leq n$. For uniformity we denote the elementary generators of ${\EQ}(2n,R,\Lambda)$ by $\eta_{ij}(*)$.

It is clear that the stabilization map takes generators of ${\EQ}(2n,R,\Lambda)$ to the
generators of ${\EQ}(2(n + 1),R,\Lambda)$.

\begin{remark}
		Throughout this article we shall assume that all ideals of $R$ are involution-invariant, i.e., if $J$ is an ideal, then $\ol{J}=J$. We also assume that $2n\geq 6$. 
\end{remark}

\begin{definition}
		The relative general quadratic subgroup of ${\GQ}(2n, R,\Lambda)$ with respect to the ideal $J$ is defined by the set $\{\alpha\in {\GQ}(2n,R,\Lambda)\mid \alpha\equiv {\rm I}_{2n} \pmod J\}$ and it is denoted by ${\GQ}(2n,R,\Lambda, J)$.
\end{definition}

\begin{definition}
		Let $(R,\Lambda)$ be a form ring and $J\subset R$ be an ideal. The subgroup of ${\GQ}(2n,R,\Lambda)$ generated by the 
		matrices of the form  $\eta_{ij}(x) \eta_{ji}(a)\eta_{ij}(x)^{-1}$, where $x\in R$ and $a\in J$, is called 
		the relative elementary subgroup and is denoted by ${\EQ}(2n,R,\Lambda, J)$.
\end{definition}

\begin{nt}
		A row $(a_1,a_2,\dots, a_n)\in R^n$ is said to be unimodular if there exists a vector $(b_1,b_2,\dots,b_n)\in\nolinebreak R^n$ such that $\sum_{i=1}^{n}a_ib_i=1$. The set of all unimodular rows of length $n$ is denoted by ${\rm Um}_n(R)$. We denote the set of all unimodular rows of length $n$ which are congruent to $e_1=(1,0,\dots,0)$ modulo the ideal $J$ by ${\rm Um}_n(R,J)$. For  an ideal $J\subset R$ the extended ideal  $J\otimes_R R[X]$ of  $R[X]$ is denoted by $J[X]$. We will mostly use localizations with respect to two types of multiplicatively
		closed subsets of $R$, {\it viz.} $S=\{1,s,s^2,\dots\}$, where $s\in R$ is a non-nilpotent
		non-zero divisor; and $S=R\setminus\mathfrak{m}$ for some $\mathfrak{m}\in \Max(R)$. By $J_s[X]$ and $J_{\mathfrak{m}}[X]$ we
		shall mean the extension of $J[X]$ in $R_s[X]$ and $R_{\mathfrak{m}}[X]$ respectively.
\end{nt}

\begin{definition}
	\label{excision}
	
		Let $(R,\Lambda)$ be a form ring and $J\subset R$ be an ideal of $R$. The excision ring of $R$ with respect to the ideal 
		$J$ is denoted by $R\oplus J$ and is defined by the set $\{(r,i)\mid r\in R,i\in J \}$ with the addition defined by $(r,i)+(s,j)=(r+s,i+j)$, and the multiplication defined by $(r,i)(s,j)=(rs,rj+is+ij)$. We can extend the involution $-$ to the ring $R\oplus J$ by setting $\ol{(r,i)}=(\ol{r},\ol{i})$. The element $(\lambda,0)\in R\oplus J$ satisfies $(\lambda,0)(\ol{\lambda},0)=(1,0)$, where $(1,0)$ is the identity element of $R\oplus J$. We can observe that  the additive subgroup $(\Lambda\oplus J)\cap \Lambda_{max}(R\oplus J)$ satisfies the properties of $(\lambda,0)$-parameter on $R\oplus J$. We fix the notation $\Gamma \oplus J$ for the subgroup $(\Lambda\oplus J)\cap \Lambda_{max}(R\oplus J)$. Hence we get the form ring $(R\oplus J,\Gamma \oplus J)$.
	
\end{definition}

Recall that there is a natural map $f: R\oplus J \rightarrow R$ given  by $f(r,i)=r+i$. This map induces a canonical homomorphism on ${\GQ}(2n,R\oplus J,\Gamma\oplus J)$. We shall use the same notation $f$ to denote this map.

The following two key lemmas are proved in \cite{BRK} for the traditional classical groups. Proofs are similar for the general 
quadratic groups.

\begin{lemma}
	\label{excision elementary}
	For $\alpha \in {\EQ}(2n, R,\Lambda,J)$, there exists a matrix $\widetilde{\alpha}\in {\EQ}(2n, R\oplus J, \Gamma \oplus J)$ such that $f(\widetilde{\alpha})=\alpha$.
\end{lemma}
%
%
\begin{lemma}
	\label{excision quadratic}
	For $\alpha\in {\GQ}(2n,R,\Lambda,J)$, there exists a matrix $\widetilde{\alpha}\in {\GQ}(2n,R\oplus J,\Gamma\oplus J)$ such that $f(\widetilde{\alpha})= \alpha$.
\end{lemma}
%
%
\section{Relative L-G Principle for the Transvection Subgroups}
It is known that any module finite ring (i.e., finite over its center) $R$ can be written as a direct limit of its  finitely generated subrings. Also, ${\G}(R, \lambda) = \underset{\lra}\lim \, {\G}(R_i, \lambda_i)$, 
where the limit is taken over all finitely generated subring of $R$. 
Hence, we can assume that $C(R)$ is Noetherian. 
For the rest of this section, we shall consider $R$ to be a module finite ring with identity.

The local-global principle for the transvection subgroups (absolute cases) of the full automorphism groups was established in \cite{BBR} for the traditional classical groups. Then, in \cite{basu2} it was generalized for the general quadratic groups. 
In this section we deduce the  relative L-G principle for the transvection subgroups.

The following results are proved in \cite{BRK} for traditional classical groups of free modules, and the steps of the proof for the general quadratic groups are identical. Therefore, we state these results without proof. 

For any column vector $v\in (R^{2n})^t$ we consider the row vector $\widetilde{v}=\ol{v}^t\psi_n$.

\begin{definition}
		We define a map $M: (R^{2n})^t\times (R^{2n})^t\rightarrow M(2n,R)$ and the inner product $\langle \:,\: \rangle$ as follows:
	$$	M(v,w) = v.\widetilde{w}-\ol{\lambda}\ol{w}\widetilde{v},~ 
	\langle v,w\rangle = \widetilde{v}.w $$

\end{definition}

\begin{lemma}
	\label{key lemma}
	Let $(R,\Lambda)$ be a form ring and $v\in {\EQ}(2n,R,\Lambda, J)e_1$. Let $w\in J^{2n}$ be
	a column vector such that $\langle v,w\rangle=0$. Then ${\rm I}_{2n}+M(v,w)\in {\EQ}(2n, R,\Lambda,J)$.
\end{lemma}

\begin{theorem}[Relative L-G principle]
	\label{lgforrelative}
		 Let $R$ be a ring and $J\subset R$ be an ideal of $R$. Let $\alpha(X)\in {\GQ}(2n, R[X], \Lambda[X], J[X])$, with $\alpha(0)={\rm I}_{2n}$ be such that for every  maximal ideal $\mathfrak{m}\in {\rm Max}(C(R))$, we have $\alpha_{\mathfrak{m}}(X)\in {\EQ}(2n, R_{\mathfrak{m}}[X], \Lambda_{\mathfrak m}[X],J_{\mathfrak{m}}[X])$. In that case,
		$\alpha(X)\in\nolinebreak {\EQ}(2n, R[X], \Lambda[X],J[X])$. 
\end{theorem}

We recall some definitions and fix notations.

\begin{definition}
		Let $(R,\Lambda)$ be a form ring and $P$ be a right $R$-module. A map $f:P\times P\to R$ is said to be a sesquilinear form if $f(pa,qb)=\ol{a} f(p,q) b$ for all  $p,q\in P$ and $a,b\in R$. A map $q:P\to R/\Lambda$ is said to be a quadratic form if $q(p)=f(p,p)+\Lambda$, where $f$ is a sesquilinear form on $P$. With respect to a sesquilinear form on $P$, we can define an associated \mbox{$\lambda$-Hermitian} form $h: P\times P\to R$ by $h(p,q)= f(p,q)+\lambda \ol{f(q,p)}$.  The triplet $(P,h,q)$ is called a quadratic module.
\end{definition}

\begin{definition}
		Let $(P,q,h)$ be a quadratic module and ${\GL}(P)$ be the full automorphism group of $P$. 
		The quadratic module $P$ is said to be non-singular if $P$ is a projective \mbox{$R$-module} and the associated $\lambda$-Hermitian form is non-singular. 
		For a non-singular $P$, the general quadratic group of $P$ is defined as follows:
		$${\GQ}(P,\Lambda,q,h)= \{\alpha\in {\GL}(P)\mid h(\alpha u,\alpha v)=h(u,v), q(\alpha u)=q(u)\}$$
		We will denote by ${\GQ}(P,\Lambda, J)$ the set $\{\alpha\in {\GQ}(P,\Lambda)\mid \alpha\equiv {\rm Id} \pmod{JP}\}$.
\end{definition}

\begin{definition}
		Let $(P,h,q)$ be  a quadratic module over $(R,\Lambda)$ and $J\subset R$ be an ideal of $R$. Let $u,v\in P$ and $a\in R$ be such that $f(u,u)\in \Lambda$, $h(u,v)=0$ and $f(v,v)=a\pmod{\Lambda}$. Then the transvection map 
		$\sigma=\sigma_{u,v,a}:P\rightarrow P$ is defined by
	$$\sigma(x)= x+uh(v,x)-v\ol{\lambda}h(u,x)-u\ol{\lambda}ah(u,x).$$
	The set of all transvections of $P$ will be denoted as ${\T}(P,\Lambda)$. A map $\sigma\in {\T}(P,\Lambda)$ is said to be a transvection relative to $J$ if either $u$ or $v$ belongs to the submodule $JP$. The set of all transvections relative to the ideal 
		$J$ will be denoted by ${\T}(P,\Lambda, J)$.
\end{definition}

\begin{definition}
		Let $(P,h,q)$ be a quadratic module over a form ring $(R,\Lambda)$ and $J\subset R$ be an ideal of $R$. Let $Q$ be the quadratic module $P\perp \mathbb{H}(R)$, where $\mathbb{H}(R)$ denotes the hyperbolic form $R\perp R^*$. Then the transvections which are of the form
	$$q=(p,a,b)\mapsto (p-aq,a,b+h(p,q)),$$
	or,
	$$q=(p,a,b)\mapsto (p-bq,a+h(p,q),b),$$
	where $a\in R$, $b\in R^*, p,q\in P$, are called elementary transvections. The set of all elementary transvections is denoted by ${\ET}(Q,\Lambda)$. An elementary transvection is said to be elementary transvection relative to $J$ if $q\in JQ$. The subgroup of ${\ET}(Q,\Lambda)$ generated by elementary transvections relative to $J$ is denoted by ${\ET}(JQ)$. And we denote ${\ET}(Q,\Lambda, JQ)$ by the normal closure of ${\ET}(JQ)$ in ${\ET}(Q,\Lambda)$. We also  use the notation ${\ET}(Q,\Lambda,J)$ for the group ${\ET}(Q,\Lambda, JQ)$.
	
\end{definition}

\begin{nt}
		\label{transvection notation}
		Let $P$ be a finitely generated quadratic $R$ module of rank $2n$ with a fixed form $\langle\:,\:\rangle$. 
		Denote by $Q$ the module $P\perp \mathbb{H}(R)$and by $Q[X]$ the module $(P\perp \mathbb{H}(R))[X]$. We assume that
		the rank of the quadratic module is $2n\geq 6$.
\end{nt}

We shall also assume the following two hypotheses:

$(H1)$ For every maximal ideal $\mathfrak{m}$ of $R$, the quadratic module $Q_{\mathfrak{m}}$ is isomorphic to $R_{\mathfrak{m}}^{2n+2}$ for the standard bilinear form $\mathbb{H}(R_{\mathfrak{m}}^{n+1})$.  

$(H2)$ For every non-nilpotent $s\in R$, if the projective module $M_s$ is a free $R_s$-module, then the quadratic module $M_s$ is isomorphic to $R_s^{2n+2}$ for the standard bilinear form $\mathbb{H}(R^{n+1})$.

\begin{lemma}[{\cite[Lemma 5.7]{RB}}, {\cite[Lemma 3.6]{BRK}}]
	\label{quadratic plus rel elementary}
	The group $${\GQ}(2n,R[X],\Lambda[X],(X))\cap {\EQ}(2n,R[X],\Lambda[X],J[X])$$ is generated by the
	elements of the form 
 \[\varepsilon\eta_{ij}(Xh(X))\varepsilon^{-1},\] 
 where $\varepsilon\in {\EQ}(2n,R,\Lambda)$,
	$h(X)\in\nolinebreak J[X]$ and $\eta_{ij}(*)$ are the elementary generators of the group ${\EQ}(2n, R[X],\Lambda[X])$.
\end{lemma}

\begin{lemma}[{\cite[Corollary 3.8]{BRK}}] 
	\label{key-dilation}
	If $\eta=\eta_1\eta_2\dots \eta_r$, where each $\eta_j$ is an elementary generator,
	and $h(Y)\in J[Y]$, then there are elements $h_t(X,Y)\in J[X,Y]$ such that $$\eta \eta_{pq}(X^{2^rm}h(Y))\eta^{-1}= \prod_{t=1}^{k}\eta_{p_tq_t}(X^mh(X,Y)).$$
\end{lemma}

We now recall some standard results:

\begin{lemma}
	\label{commutes with poly}
	Let $R$ be a ring and $K$ a finitely presented left (right) $R$-module, and let $L$ be any  left (right) $R$-module. 
	Then we have a natural isomorphism:
	$$f: {\rm Hom}_R(K,L)[X] \rightarrow {\rm Hom}_{R[X]}(K[X],L[X]).$$
\end{lemma}

\begin{lemma}
	\label{commutes with Hom}
	Let $S$ be a multiplicatively close subset of a ring $R$. Let $K$ be a finitely presented $R$-module and $L$ be any $R$-module. Then we have a natural isomorphism $$g: S^{-1}({\rm Hom}_R(K,L)) \rightarrow {\rm Hom}_{S^{-1}R}(S^{-1}K,S^{-1}L).$$
\end{lemma}

The following lemma is used frequently (sometimes in a subtle way) in the proof of the main results. 

\begin{lemma}[{cf. \!\cite[Lemma 5.1]{HV}}]
	\label{noeth} 
	Let $A$ be Noetherian ring and $0\neq s\in A$. 
	Then there exists a natural number 
	$k$ such that the homomorphism \[{\G}(A,s^kA, s^k\lambda) \ra {\G}(A_s, \lambda_s)\] 
	induced by the localization homomorphism $A \ra A_s$ is injective.  
\end{lemma}

Now we prove the relative L-G principle by using Lemma \ref{excision elementary}. 

\begin{proposition}[Relative Dilation Principle]
	\label{rel-dilation-trans}
	Let $R$ be an almost commutative ring (i.e., an 
	associative ring which is finite over its center $C(R)$) and $J\subset R$ be an ideal. Let $P$ and $Q$ be as in \ref{transvection notation}. Let $s$ be a non-nilpotent element of $R$ such that $P_s$ is a free module.  Let $\sigma(X)\in {\GQ}(Q[X], \Lambda[X],J[X])$ with $\sigma(0)=Id$and suppose additionally that $\sigma_s(X)\in {\EQ}(2n+2,R_s[X],\Lambda_s[X],J_s[X])$. Then there exists $\widehat{\sigma}(X)\in {\ET}(Q[X],\Lambda[X], J[X])$ and $l>0$ such that $\widehat{\sigma}(X)$ localizes to $\sigma(bX)$ for some $b\in (s^l)$ and $\widehat{\sigma}(0)=Id$.
\end{proposition}

\begin{proof} Since elementary transvections can always be lifted, then we may assume that $R$ is reduced. We will show that there exists $l>0$ such that $$\sigma(bX)\in {\ET}(Q[X],\Lambda [X], J[X])$$ for all $b\in (s^k)$ and for all $k\geq l$.

As $\sigma(0)= {\rm Id}$, by Lemma \ref{quadratic plus rel elementary}, we can write $\sigma_s(X)= \prod_{k}\gamma_k\eta_{i_kj_k}(X\lambda_k(X))\gamma_k^{-1}$, where $\gamma_k\in {\EQ}(2n+2,R_s,\Lambda_s)$, and $\lambda_k(X)\in J_s[X]$. Hence, by  the proof of Lemma \ref{excision elementary}, there exists $\widetilde{\sigma}_{(s,0)}(X)\in {\EQ}(2n+2, (R_s\oplus J_s)[X],(\Gamma_s\oplus J_s)[X])$ such that
$$\widetilde{\sigma}_{(s,0)}(X)= \prod_{k}\widetilde{\gamma_k}\widetilde{\eta_{i_kj_k}}(0,X\lambda_k(X))\widetilde{\gamma_k}^{-1},$$
where $\phi_s(\widetilde{\gamma_k})=\gamma_k$, $(0,X\lambda_k(X))\in (R\oplus J)_{(s,0)}[X]$, the $\widetilde{\eta_{i_kj_k}}(*)$ are the elementary generators of ${\EQ}(2n+2, (R_s\oplus J_s)[X],(\Gamma_s\oplus J_s)[X])$  and $\phi:R\oplus J \rightarrow R$ is defined by $\phi((a,i))=a+i$ and $\Gamma_s\oplus J_s=(\Gamma \oplus J)_s$ is the localization of $\Gamma\oplus J$ with respect to the set $\{1,s,s^2,\dots\}$. Hence, for $d>0$, we have that 
$$\widetilde{\sigma}_{(s,0)}(XT^{2d})= \prod_{k}\widetilde{\gamma}_k\widetilde{\eta_{i_kj_k}}((0,XT^{2d}\lambda_k(XT^{2d}))\widetilde{\gamma}_k^{-1},$$
for some $\gamma_k\in {\EQ}(2n+2, R_s,\Lambda_s)$. Using Lemma \ref{key-dilation} and standard commutator formulas (see \cite[Lemma 3.16, pg. 43]{Bak1}), we get that $\widetilde{\sigma}_{(s,0)}(XT^{2d})= \prod_{t}\widetilde{\eta_{p_tq_t}}(T\mu_t(X))$, for some $\mu_t(X)\in (R\oplus J)_{(s,0)}[X]$ with $p_t=1$ or $q_t=1$.

Since $P_s$ is a free $R_s$-module, then we have 
$$(P\oplus J)_{(s,0)}[X,T] \cong (R\oplus J)_{(s,0)}^{2n}[X,T] \cong (P\oplus J)_{(s,0)}[X,T]^*.$$
Thus using the isomorphism, polynomials in $(P\oplus J)_{(s,0)}[X,T]$ can be regarded as linear forms. 

First we consider the case: $p_t=1$. Let $p_1^*,p_2^*,\dots, p_n^*,p_{-1}^*,\dots, p_{-n}^*$ be the standard basis of $(P\oplus J)_{(s,0)}$. Let $s^mp_i^*\in P\oplus J$ for some $m>0$ and $i=\pm 1,\pm 2,\dots,\pm n$. Let $e_{\pm i}^*$ be the standard basis of $(R\oplus J)^{2n}$. Then for $q_t=\pm i$, consider the element $T\mu_t(X)e_{\pm i}^*$ as an element in $(P\oplus I)_{(s,0)}[X,T]^*$. As $(P\oplus J)_{(s,0)}$ is free, by Lemma \ref{commutes with poly}, we can say $T\mu_t(X)e_{\pm i}^*$ is a polynomial in $T$. Again, by Lemma \ref{commutes with Hom}, there exists $k_1>0$ such that $k_1$ is the maximum power of $(s,0)$ occurring  in the denominator of $\mu_t(X)e_{\pm i}^*$. Choose $l_1\geq \max(k,m)$.

Now consider the case $q_t=1$. Then, for $p_t=\pm j$, $T\mu_t(X)e_{\pm j}^*\in (P\oplus J)_{(s,0)}[X,T]$. By a similar argument we can consider $T\mu_t(X)e_{\pm j}^*$ as a polynomial in $T$ and hence there exists $k_2>0$ such that $k_2$ is maximum power of $(s,0)$ occurred in $\mu_te_{\pm j}^*$. Now choose $l_2\geq\nolinebreak \max(k_2,m)$. For $l\geq \max(l_1,l_2)$, under the transformation $T\mapsto (s,0)^lT$, $\widetilde{\sigma}_{(s,0)}((b,0)XT^{2d})$ is defined  over $(Q\oplus J)[X,T]$, i.e., there exists some element 
$$\widetilde{\widehat{\sigma}}(X,T)\in {\ET}((Q\oplus J)[X,T])$$
such that $\widetilde{\widehat{\sigma}}_{(s,0)}(X,T)= \widetilde{\sigma}_{(s,0)}((b,0)XT^{2d})$. Replacing $T=(1,0)$ and using Lemma~\ref{noeth}, one gets $\widetilde{\sigma}((b,0)X)\in {\ET}(Q[X]\oplus J[X],(\Gamma \oplus J)[X], 0\oplus J[X])$. Hence, the result follows by applying $\phi$. \end{proof}

Consequence: {\bf Relative L-G principle for the transvection subgroups}: 

\begin{theorem}
	\label{lg for rel trans}
	Let $R$ be an almost commutative ring and $J\subset R$ be an ideal. Let $P$ and $Q$ be as in \ref{transvection notation}. Let 
	$\sigma(X)\in {\GQ}(Q[X],\Lambda[X], J[X])$ with $\sigma(0)=Id$.
 
 If $\sigma_{\mathfrak{m}}(X) \in {\EQ}(2n+2, R_{\mathfrak{m}}[X],\Lambda_{\mathfrak{m}}[X], J_{\mathfrak{m}}[X])$ for all $\mathfrak{m}\in {\rm Max}(C(R))$, then we have $\sigma(X)\in {\ET}(Q[X],\Lambda[X],J[X])$. 
\end{theorem}

\begin{proof} Follows by arguing as in the proof of \cite[Lemma 3.10]{basu2}, and using Proposition \ref{rel-dilation-trans}.
\end{proof}

\section{Relative Stability for Quadratic ${\kk}$}

The aim of this section is to establish the ${\kk}$-stability of the relative transvection groups as an application of Theorem \ref{lg for rel trans}. For the absolute case we refer to \cite{BTP} and \cite{basu2}. 

\begin{definition}
		\label{double ring}
		Let $R$ be an associative ring with identity and $J\subset R$ an ideal. Consider the ring $D$ = $\{(a,b)\in R\times R: a-b\in J\}$ with addition and multiplication defined component wise. We call it the double ring relative to the ideal $J$. For a form ring $(R,\Lambda)$, one extends the involution of $-: R\rightarrow R$ to the ring $D$, defining $-:D\rightarrow D$ by $\ol{(a,b)}=(\ol{a},\ol{b})$. We fix the element $(\lambda,\lambda)$ and define $\Lambda^{\prime}$ = $\{(a,b)\in \Lambda \times \Lambda\mid a-b\in J \}$. Then one can show that $(D,\Lambda^{\prime})$ is a form ring.
\end{definition}

\begin{definition}
		Let $(R,\Lambda,\lambda)$ and $(S,\Lambda^{\prime},\lambda^{\prime})$ be two form rings.  A ring homomorphism $f:R\rightarrow S$ is said to be a morphism of form rings if $f(\ol{r})=\ol{f(r)}$, $f(\Lambda)\subseteq \Lambda^{\prime}$ and $f(\lambda)=\lambda^{\prime}$.
\end{definition}

\begin{lemma}
	\label{isomorphic excision}
	Let $(R,\Lambda)$ be a form ring and $J\subset R$ be a two sided  ideal. Then the form rings $(D,\Lambda^{\prime})$ and $(R\oplus J, \Lambda \oplus J)$ are isomorphic. 
\end{lemma}

\begin{proof} Consider the homomorphisms $f: D\rightarrow R\oplus J$ defined by $f(a,b)=(a,b-a)$ and $g:R\oplus J\rightarrow D$ defined by $g(a,i)=(a,a+i)$. It can be checked that both $f$ and $g$ are form homomorphisms. They are inverses of each other. \end{proof}

\begin{lemma}[{\cite[Proposition 3.1]{keshari}}]
	\label{dimension of excision}
	Let $A$ be a commutative Noetherian ring of (Krull) dimension $d$ and $J\subset A$ be an ideal. Then the ring $D$  is also a commutative ring of dimension $d$.
\end{lemma}

\begin{proof} Clearly $D_A$ is a commutative Noetherian ring. We first prove that dimension of $A\oplus J$ is $d$, and the rest follows from Lemma \ref{isomorphic excision}. Clearly $A$ can be identified with the sub-ring $\{(r,0)\mid r\in A\}$ of $A\oplus J$. The element $(0,i)$ is integral over $A$ since we have that $(0,i)^2-(i,0)(0,i)=(0,0)$. Hence every element of $A\oplus J$ is integral over $R$, and therefore $\dim(A)=\dim(A\oplus J)$. \end{proof} 

Recall ${\kk}$-stability result of Bak--Petrov--Tang (cf. \!\cite{BTP})
for general quadratic groups in the absolute case.	
\begin{theorem}[{\cite{BTP}}]
	\label{Global stability}
	
	Let $R$ be an almost commutative with $\dim(C(R))=d$. Consider the form ring $(R, \Lambda)$. 
	Then the stabilization map 
	$$\frac{{\GQ}(2n, R,\Lambda)}{{\EQ}(2n, R,\Lambda)}\rightarrow \frac{{\GQ}(2n+2, R, \Lambda)}{{\EQ}(2n+2, R, \Lambda)}$$ 
	is an isomorphism for $2n\geq \max(6,2d+4)$.
\end{theorem}

We prove the above result in the  relative case.

\begin{theorem}
	\label{relative stability}
	Let $R$ be a form ring which is finitely generated  over its center $C(R)$, with $\dim(C(R))=d$ and let $J\subset R$ be a two-sided ideal of $R$. Then the stabilization map 
	$$\frac{{\GQ}(2n, R,\Lambda,J)}{{\EQ}(2n, R,\Lambda,J)}\rightarrow \frac{{\GQ}(2n+2, R,\Lambda,J)}{{\EQ}(2n+2, R, \Lambda,J)}$$ is an isomorphism for $2n\geq \max(6,2d+4)$.
\end{theorem}

\begin{proof}  
Consider the stabilization map
$\phi: {\KQ}_{1,2n}\rightarrow {\KQ}_{1,2n+2}$, where $${\KQ}_{1,2n}(R,\Lambda,J)=\frac{{\GQ}(2n, R,\Lambda,J)}{{\EQ}(2n,R,\Lambda,J)}.$$

By \cite[5.3.22]{HO}, we have following sequences are exact
\[
\begin{tikzcd}
	&1\arrow{r}{} &{\GQ}(2n,R,\Lambda,J)\arrow{r}{i}&{\GQ}(2n,D,\Lambda^{\prime})\arrow{r}{p_2} &{\GQ}(2n, R,\Lambda) \arrow{r}{}&1, \\
	&1\arrow{r}{} &{\EQ}(2n,R,\Lambda,J)\arrow{r}{i}&{\EQ}(2n,D,\Lambda^{\prime})\arrow{r}{p_2} &{\EQ}(2n, R,\Lambda) \arrow{r}{}&1,
\end{tikzcd}
\]
where the map $i$ is induced from the map $M_{2n}(J)\rightarrow M_{2n}(D)$ given by $i(\alpha)=(\alpha,{\rm I}_{2n})$ and the map $p_2$ is induced from the projection map $p_2: D\rightarrow R$ given by $p_2((a,b))=b$. Thus we have the following commutative diagram:

\[
\begin{tikzcd}
	& &1\arrow{d}{} &1\arrow{d}{} &1\arrow{d}{} & \\
	&1\arrow{r}{} &{\EQ}(2n,R,\Lambda,J)\arrow{r}{}\arrow{d}{} &{\EQ}(2n,D,\Lambda^{\prime})\arrow{r}{} \arrow{d}{} &{\EQ}(2n, R,\Lambda) \arrow{r}{} \arrow{d}{} &1 \\
	&1\arrow{r}{}&{\GQ}(2n, R,\Lambda,J)\arrow{r}{} \arrow{d}{} &{\GQ}(2n, D,\Lambda^{\prime})\arrow{r}{} \arrow{d}{}&{\GQ}(2n, R,\Lambda)\arrow{r}{} \arrow{d}{}&1\\
	&1\arrow{r}{} &{\KQ}_{1,2n}(R,\Lambda,J)\arrow{r}{} \arrow{d}{} &{\AQ}_{1,2n}(D,\Lambda^{\prime}) \arrow{r}{} \arrow{d}&{\AQ}_{1,2n}(R,\Lambda) \arrow{r}{}\arrow{d}{} &1\\
	&  &1 &1 &1
\end{tikzcd}
\]

Taking limits over $n$ in the lower row, we get an exact sequence 

\[\begin{tikzcd}
	&1\arrow{r}{} &{\KQ}_1(R,\Lambda,J)\arrow{r}{}&{\KQ}_1(D,\Lambda^{\prime})\arrow{r}{} &{\KQ}_1(R,\Lambda) \arrow{r}{}&1,
\end{tikzcd}\]

and for any $n$, we get a homomorphism of exact sequences:

\[
\begin{tikzcd}
	&1\arrow{r}&{\KQ}_{1,2n}(R,\Lambda,J)\arrow{r}\arrow{d} &{\KQ}_{1,2n}(D,\Lambda^{\prime}) \arrow{r}\arrow{d}
	&{\KQ}_{1,2n+2}( R,\Lambda)\arrow{d}\arrow{d}\arrow{r} &1\\
	&1\arrow{r} &{\KQ}_1( R,\Lambda,J) \arrow{r} &{\KQ}_1(D,\Lambda^{\prime})\arrow{r}&{\KQ}_1(R,\Lambda)\arrow{r}&1.
\end{tikzcd}
\]
For $2n\geq 2d+4$ the right most map is an isomorphism by Theorem \ref{Global stability}. Since $R$ is finitely generated over $C(R)$, then it can be checked that $D$ is finitely generated over $C(D)$. Since $C(D)$ is a double ring of $C(R)$ relative to the ideal 
$C(R)\cap J$, by Lemma \ref{dimension of excision}, we get $\dim(C(D))=d$. Hence for $2n\geq 2d+4$, the middle homomorphism is also an isomorphism by Theorem \ref{Global stability}. Hence it follows that the left most map is an isomorphism. \end{proof}

Now we prove the ${\kk}$-stability of relative transvection groups as an application of relative local-global principle of transvection groups. For this we need to recall the following result of Vaserstein.

\begin{lemma}[{\cite[Lemma 4.2]{basu2}}]
	\label{key-stability}
	Let $(R,\Lambda)$ be a form ring finitely generated over $C(R)$ with Krull dimension of $C(R)=d$ and $J\subset R$ be an  two-sided ideal of $R$. Let $P,Q$ be as in Notation \ref{transvection notation}. Let the rank of $Q$ be $2n\geq {\rm Max}(6, 2d+2)$. Then  the group of elementary transvections ${\ET}(Q\perp \mathbb{H}(R),J)$ acts transitively on the set ${\rm Um}(Q\perp \mathbb{H}(R),J)$ of unimodular elements which are congruent to $(0,\dots,0,1,0)$ modulo $J$.   
\end{lemma}

\begin{theorem}
	Let $(R,\Lambda)$ be a form ring finitely generated over $(C(R))$ with Krull dimension $d$ and let $J\subset R$ be a two-sided ideal. Let $P,Q$ be as defined in Notation \ref{transvection notation} and let the rank of $Q$ be $2n\geq \max(6,2d+4)$. Then the stabilization map $$ i_{2n}:{\AQ}_{1,2n}(Q,J) \rightarrow {\AQ}_{1,2n+2}(Q\perp \mathbb{H}(R),J)$$ is an isomorphism.  
\end{theorem}

\begin{proof} 
The surjectivity part follows from Lemma \ref{key-stability}. In order to prove the injectivity part, let $\alpha\in {\GQ}(Q,J)$ be such that $\widetilde{\alpha}=\alpha\perp \rm Id$ lies in ${\ET}(M\perp \mathbb{H}(R),J)$. Let $\varphi(X)$ be the isotropy between $\widetilde{\alpha}$ and $\rm Id$. Now by similar argument as given in \cite[Theorem 4.4]{basu2}, we can get an isotropy $\widetilde{\varphi}(X)\in {\GQ}(Q[X],J[X])$ between $\alpha$ and $\rm Id$. Now, by localizing $\widetilde{\varphi}$ 
at a maximal ideal $\mathfrak{m}\in \Max(C(R))$, $\widetilde{\varphi}_{\mathfrak{m}}(X)$ becomes (relative) stably elementary. 
By Theorem~\ref{relative stability} it follows that 
$\widetilde{\varphi}_{\mathfrak{m}}(X)$ is actually (relative)  elementary. 
Therefore, by the relative L-G principle (Theorem \ref{lg for rel trans}) one gets $\widetilde{\varphi}(X)\in {\ET}(Q[X]),J[X])$. Hence 
$\alpha=\nolinebreak\widetilde{\varphi}(1)\in\nolinebreak {\ET}(Q,J)$. \end{proof}

\section{Structure of Unstable Quadratic ${\kk}$ Groups; Relative Case}

We devote this section to discuss the study of nilpotent property of (relative) unstable
${\rm K_1}$-groups. Throughout this section we assume $R$ is a commutative ring with identity, i.e., we are considering the trivial involution and $n\ge 3$. By ${\SQ}(2n, R,\lambda)$ we shall denote the subgroup of ${\GQ}(2n, R,\lambda)$ with matrices of determinant 1, and analogously for the relative cases. 

Let $G$ be a group. Define $Z^i$ inductively by $Z^0=G$, $Z^1=[G,G]$ and $Z^i=[G,G^{i-1}]$. The group $G$ is said to be nilpotent if $Z^r=\{e\}$ for some $r>0$.

\begin{definition}
	A group $G$ is called nilpotent-by-abelian if it has a normal subgroup $H$ such that $H$ is nilpotent and $G/H$ is abelian. 
\end{definition}

In  \cite{Bak}, A.Bak proved that  
the unstable ${\kk}(R)$ group of ${\GL}_n(R)$ is nilpotent-by-abelian for $n\ge 3$, and hence ${\kk}(R)$ is solvable. 
Later it was generalized by R. Hazrat for the general quadratic groups over module finite rings (cf. \!\cite{HR}).
In \cite{HV} Hazrat--Vavilov revisited this problem for ordinary classical Chevalley groups (that is, types A, C, and D)
and finally extends it further to the exceptional Chevalley groups (that is, types E, F, and G). 
A simpler and shorter proof is given 
in \cite{BBR} for the linear, symplectic and orthogonal groups (absolute cases). In \cite{RB}, 
the first author proved this result for the general Hermitian groups and the same proof also works for the general quadratic groups. The relative cases are proved in \cite{BBR} for the traditional classical groups. In this article we 
prove the result for the relative general quadratic groups, and consequently we get the result for the module cases as an 
application of (relative) L-G principle for the transvection subgroups.

\begin{theorem}
	\label{relative nilpotent}
	The quotient group $\SQ(2n, R,\Lambda,J)/\EQ(2n, R,\Lambda,J)$ is nilpotent for $n\geq 3$. The class of nilpotency is at most $\max(1, d+2-n)$ where $d=\dim R.$ 
\end{theorem}

\begin{proof} 
Note that for $n\geq d+2$, the quotient group

$$G={\SQ}(2n, R,\Lambda,J)/{\EQ}(2n, R,\Lambda,J)$$
is abelian and hence nilpotent. So we consider the case $n< d+2$. Let us fix a natural number $n$. We prove the theorem by induction on $d=\dim R$. 
Let  $m= d+2-n$ and $\alpha=[\beta, \gamma]$ for some $\beta \in G$ and $\gamma \in Z^{m-1}$. Clearly the result is true for $d=0$. Let $\widetilde{\beta}$ be the pre-image of $\beta$  and $\widetilde{\gamma}$ 
be the pre-image of $\gamma$ under the map $${\SQ}(2n, R,\Lambda,J)\rightarrow G.$$

Consider the double ring $(D,\Lambda^{\prime})$ as defined in \ref{double ring}. By \cite[5.3.22]{HO}, we have the following exact sequences:
\[
\begin{tikzcd}
	&1\arrow{r}{} &{\SQ}(2n,R,J,\Lambda)\arrow{r}{i}&{\SQ}(2n,D,\Lambda^{\prime})\arrow{r}{p_2} &{\SQ}(2n, R,\Lambda) \arrow{r}{}&1 \\
	&1\arrow{r}{} &{\EQ}(2n,R,J,\Lambda)\arrow{r}{i}&{\EQ}(2n,D,\Lambda^{\prime})\arrow{r}{p_2} &{\EQ}(2n, R,\Lambda) \arrow{r}{}&1
\end{tikzcd}
\]
where the map $i$ is induced from the map $M_{2n}(J)\rightarrow M_{2n}(D)$ given by $i(\alpha)=(\alpha,\rm I_{2n})$ and $p_2$ is induced from the projection map $p_2: D\rightarrow R$ given by $p_2((a,b))=b$. Now, in the group $H={\SQ}(2n, D,\Lambda^{\prime})/{\EQ}(2n, D,\Lambda^{\prime})$ we have $i(\widetilde{\gamma})\in H^m$. Since ${\rm dim}(R)={\rm dim(D)}$ (by Lemma \ref{dimension of excision}), by the result of the absolute cases, we get $[i(\widetilde{\beta}),i(\widetilde{\gamma})]\in {\EQ}(2n, D, \Lambda^{\prime})$. It can be checked that $i([\widetilde{\beta},\widetilde{\gamma}])=[i(\widetilde{\beta}),i(\widetilde{\gamma})]$ and the image of this element under $p_2$ is identity. Hence by the above diagram one gets $[\widetilde{\beta},\widetilde{\gamma}]\in {\EQ}(2n, R,J)$, 
consequently $[\beta,\gamma]$ is trivial in $G$. Next, consider the pre-images $\beta_1$ of $\widetilde{\beta}$ and  $\gamma_1$ of $\widetilde{\gamma}$ under the map ${\SQ}(2n,R\oplus J,\Gamma\oplus J)\rightarrow\nolinebreak {\SQ}(2n, R,\Lambda,J)$. Now by the result of absolute cases we get $[\beta_1,\gamma_1]\in\nolinebreak EQ(2n, R\oplus J,\Gamma \oplus J)$. It can also be checked that the same commutator $[\beta_1, \gamma_1]$ is in ${\EQ}(2n, R\oplus J,\Gamma \oplus J,0\oplus J)$. By projecting $R\oplus J$ onto $R$, 
one gets $[\widetilde{\beta},\widetilde{\gamma}]\in {\EQ}(2n, R,J,\Lambda)$, and 
hence $\alpha=\{1\}$ in $G$. \end{proof}

\begin{corollary}
	Let $(R,\Lambda)$ be a commutative form ring and $J\subset R$ be an ideal of $R$. Then the quotient group ${\GQ}(2n, R,\Lambda,J)/{\ET}(2n, R,\Lambda,J)$ is nilpotent-by-abelian for $n\geq 3$.
\end{corollary}
%

\begin{corollary}
	Let $(R,\Lambda)$ be a commutative form ring and $J\subset R$ be an ideal of $R$. Consider the notation as in \ref{transvection notation}. Let $d=\dim(R)$ and $t= $ the local rank of $Q$. The quotient group ${\T}(Q,J)/{\ET}(Q,J)$ is nilpotent of class at most $\max(1, d+3-t/2)$.
\end{corollary}

\begin{proof} The proof is same as \cite[Theorem 4.1]{BBR}. \end{proof}

\section{Bass' Nil Group ${\rm NKQ_1(R)}$}

The group ${\nk}(R)=\mathsf{ker}
({\kk}(R[X])\ra {\kk}(R))$; $X=0$ is called the Bass' nil-group of R. This is the subgroup consisting of elements $[\alpha(X)]\in {\kk}(R[X])$ 
such that $[\alpha(0)]=[{\I}]$, and hence ${\kk}(R[X])\cong {\rm NK}_1(R)\oplus {\kk}(R)$. 
The aim of the next sections is to study some properties of Bass nil-groups ${\rm NK_1}$ for 
the general quadratic groups or Bak's unitary groups.

In this section we recall some basic definitions and properties of the representatives 
of ${\rm NKQ_1}(R)$.
We represent any element of
${\rm M}_{2n}(R)$ as a matrix 
\[\begin{pmatrix}
		a&b\\c&d
	\end{pmatrix}, \text{where }a,b,c,d\in {\rm M}_n(R).\] 
 For a matrix as above, we call {\small $\begin{pmatrix}
		a&b
	\end{pmatrix}$} the upper half of $\alpha$. Let $(R,\lambda, \Lambda)$ be a form ring. 
By setting $\ol{\Lambda}=\{\ol{a}\mid a\in \Lambda\}$ we get another form ring $(R,\ol{\lambda},\ol{\Lambda})$.
We can extend the involution of $R$ to ${\rm M}_n(R)$ by setting $(a_{ij})^*=(\ol{a}_{ji})$. For details, see \cite{V}, \cite{V1}. 

\begin{definition}
	Let $(R,\lambda,\Lambda)$ be a form ring. A matrix $\alpha=(a_{ij})\in {\rm M}_n(R)$ 
		is said to be $\Lambda$-Hermitian if $\alpha=-\lambda \alpha^*$ and all the diagonal entries of $\alpha$ 
		are contained in  $\Lambda$. A matrix $\beta\in {\rm M}_n(R)$ is said to be $\ol{\Lambda}$-Hermitian 
		if $\beta=-\ol{\lambda}\beta^*$ and all the diagonal entries of $\beta$ are contained in $\ol{\Lambda}$. 
\end{definition}

\begin{remark}
	A matrix $\alpha\in {\rm M}_n(R)$ is $\Lambda$-Hermitian if and only if $\alpha^*$ is $\ol{\Lambda}$-Hermitian.   
\end{remark}

\begin{lemma}[{\cite[Example 2]{V}}]
	\label{stableunderconjugation}
	Let $\beta\in {\GL}_n(R)$ be a $\Lambda$-Hermitian matrix. Then the matrix $\alpha^*\beta \alpha$ is $\Lambda$-Hermitian for every $\alpha \in {\GL}_n(R)$.
\end{lemma}

\begin{definition}
	\label{lamdaunitary}
	Let $\alpha=\begin{pmatrix}
			a&b\\c&d
		\end{pmatrix}\in {\rm M}_{2n}(R)$ be a matrix. Then $\alpha$ is said to be a $\Lambda$-quadratic matrix if one of the following equivalent conditions holds:
	\begin{enumerate} 
		\item $\alpha \in {\GQ}(2n, R,\Lambda)$ and the diagonal entries of the matrices $a^*c$, $b^*d$ are in $\Lambda$,
		\item $a^*d+\lambda c^*b={\I}_n$ and the matrices $a^*c$, $b^*d$ are $\Lambda$-Hermitian, \label{twolambdaunitary}
		\item $\alpha\in {\GQ}(2n,R,\Lambda)$ and the diagonal entries of the matrices $ab^*,cd^*$ are in $\Lambda$,
		\item$ad^*+\lambda bc^*={\I}_{n}$ and the matrices $ab^*$, $cd^*$ are $\Lambda$-Hermitian. \label{fourlambdaunitary}
	\end{enumerate}
\end{definition}


\begin{lemma}
	Let $\alpha=\begin{pmatrix}
		a&0\\0&d
	\end{pmatrix}\in {\rm M}_{2n}(R)$. Then $\alpha\in {\GQ}^{\lambda}(2n, R,\Lambda)$ if and only if we have $a\in {\GL}_n(R)$ and $d=(a^*)^{-1}$.
\end{lemma}

\begin{proof} Let $\alpha \in {\GQ}^{\lambda}(2n, R, \Lambda)$. In view of condition \ref{twolambdaunitary}.\ of Definition \ref{lamdaunitary}, we have, $a^*d= {\I}_n$. Hence $a$ is invertible 
and $d=(a^*)^{-1}$. The converse holds by condition \ref{twolambdaunitary}.\ of Definition~\ref{lamdaunitary}. \end{proof}

\begin{definition}
	Let $\alpha\in {\GL}_n(R)$ be a matrix. A matrix of the form {\small $\begin{pmatrix}
				\alpha&0\\0&(\alpha^*)^{-1}
			\end{pmatrix}$} is denoted by $\mbb{H}(\alpha)$ and is said to be hyperbolic.
\end{definition}

\begin{remark}
	In a similar way we can show that a matrix of the form  $T_{12}(\beta):=\begin{pmatrix}
			{\I}_n &\beta\\0&{\I}_n
		\end{pmatrix}$ is a $\Lambda$-quadratic matrix if and only if $\beta$ is $\ol{\Lambda}$-Hermitian. Similarly, a matrix of the form $T_{21}(\gamma):=\begin{pmatrix}
			{\I}_n&0\\\gamma&{\I}_n
		\end{pmatrix}$ is a $\Lambda$-quadratic matrix if and only if $\gamma$ is $\Lambda$-Hermitian.
\end{remark}

Likewise, in the quadratic case we can define the notion of $\Lambda$-elementary quadratic groups in the following way:

\begin{definition}
		The $\Lambda$-elementary quadratic group is denoted by ${\EQ}^{\lambda}(2n, R, \Lambda)$ 
		and defined to be the group generated by the $2n\times 2n$ matrices of the following forms:
  \begin{itemize}
      \item $\mbb{H}(\alpha)$, here $\alpha\in {\E}_n(R)$,
      \item $T_{12}(\beta)$, where $\beta$ is $\ol{\Lambda}$-Hermitian, and
      \item $T_{21}(\gamma)$, where $\gamma$ is $\Lambda$-Hermitian.
  \end{itemize}  
\end{definition}

\begin{lemma}
	\label{uppertriangular}
	Let $A=\begin{pmatrix}
		\alpha&\beta\\0&\delta
	\end{pmatrix}\in {\rm M}_{2n}(R)$. Then $A\in {\GQ}^{\lambda}(2n,R,\Lambda)$ if and only if we have
	$\alpha\in {\GL}_n(R)$, $\delta=(\alpha^*)^{-1}$ and $\alpha^{-1}\beta$ is $\ol{\Lambda}$-Hermitian. In this case 
	$$A\equiv \mbb{H}(\alpha) \pmod{{\EQ}^{\lambda}(2n, R, \Lambda)}.$$ 
\end{lemma}

\begin{proof} Let $A\in {\GQ}^{\lambda}(2n, R, \Lambda)$. Then, by condition \ref{fourlambdaunitary}.\ of Definition \ref{lamdaunitary}, we have that 
$\alpha\delta^*={\I}_n$ and $\alpha \beta^*$ is $\Lambda$-Hermitian. Hence $\alpha$ is invertible and 
$\delta=(\alpha^*)^{-1}$. For $\alpha^{-1}\beta$, we get $$(\alpha^{-1}\beta)^*=\beta^*(\alpha^{-1})^*=\alpha^{-1}(\alpha \beta^*)(\alpha^{-1})^*,$$ 
which is $\Lambda$-Hermitian by Lemma \ref{stableunderconjugation}. Hence $\alpha^{-1}\beta$ is $\ol{\Lambda}$-Hermitian. 
Conversely, condition on $A$ fulfills condition \ref{fourlambdaunitary}.\ of Definition \ref{lamdaunitary}. 
Hence $A$ is $\Lambda$-quadratic. Now  
, 
since $\alpha^{-1}\beta$ is $\ol{\Lambda}$-Hermitian, $T_{12}(-\alpha^{-1} \beta)\in {\EQ}^{\lambda}(2n, R,\Lambda)$,
and $AT_{12}(\alpha^{-1} \beta)= \mbb{H}(\alpha)$. 
Thus $A\equiv \mbb{H}(\alpha) \pmod{{\EQ}^{\lambda}(2n, R,\Lambda)}$. \end{proof}  

Arguing similarly one gets the following:

\begin{lemma}
	\label{lowertriangular}
	Let $B=\begin{pmatrix}
		\alpha &0\\ \gamma &\delta
	\end{pmatrix}\in {\rm M}_{2n}(R)$. Then $B\in {\GQ}^{\lambda}(2n, R,\Lambda)$ if and only if 
	$\alpha\in {\GL}_n(R)$, $\delta=(\alpha^*)^{-1}$ and $\gamma$ is $\Lambda$-Hermitian. 
	In this case  $$B\equiv \mbb{H}(\alpha) \pmod{{\EQ}^{\lambda}(2n, R, \Lambda)}.$$ 
\end{lemma}

\begin{lemma}
	\label{ifaisinvertible}
	Let $\alpha=\begin{pmatrix}
		a&b\\
		c&d
	\end{pmatrix}\in {\GQ}^{\lambda}(2n,R,\Lambda)$. If $a\in {\GL}_n(R)$ then $$\alpha \equiv \mbb{H}(a) \pmod{{\EQ}^{\lambda}(2n,R,\Lambda)}.$$
\end{lemma}

\begin{proof} By same argument as given in Lemma \ref{uppertriangular}, we have $a^{-1}b$ is $\ol{\Lambda}$-Hermitian. 
Hence  $T_{12}(-a^{-1}b)\in {\EQ}^{\lambda}(2n, R, \Lambda)$, and consequently $\alpha T_{12}(-a^{-1}b)=\begin{pmatrix}
	a&0\\c &d^{\prime}
\end{pmatrix}\in {\GQ}^{\lambda}(2n, R, \Lambda)$ for some $d^{\prime}\in {\GL}_n(R)$. Hence by Lemma \ref{lowertriangular}, we get 
$$\alpha T_{12}(-a^{-1}b)\equiv H(a) \pmod{{\EQ}^{\lambda}(2n, R, \Lambda)}.$$ 
Hence  $\alpha\equiv H(a)\pmod{{\EQ}^{\lambda}(2n, R, \Lambda)}$. \end{proof}

\begin{definition}Consider the  embedding:
	$${\GQ}^{\lambda}(2n, R, \Lambda)\rightarrow {\GQ}^{\lambda}(2n+2,R,\Lambda), \,\,\,$$ 
	We use the following notation: 
 $${\GQ}^{\lambda}(R,\Lambda)=\underset{n=1}{\overset{\infty}\cup}{\GQ}^{\lambda}(2n, R,\Lambda),$$ 
		$${\EQ}^{\lambda}(R,\Lambda)=\underset{n=1}{\overset{\infty}\cup} {\EQ}^{\lambda}(2n, R,\Lambda).$$
\end{definition}

In view of quadratic analog of Whitehead Lemma, we have that the group ${\EQ}^{\lambda}(R,\Lambda)$ 
coincides with the commutator of ${\GQ}^{\lambda}(R, \Lambda)$. Therefore, the group 
$${\rm KQ_1}^{\lambda}(R,\Lambda):= \frac{{\GQ}^{\lambda}(R,\Lambda)}{{\EQ}^{\lambda}(R,\Lambda)}$$ is well-defined. The class of a matrix 
$\alpha\in {\GQ}^{\lambda}(R,\Lambda)$ in the group ${\rm KQ_1}^{\lambda}(R,\Lambda)$ is denoted by $[\alpha]$. In this way we obtain a 
${\kk}$-functor ${\rm KQ_1}^{\lambda}$ acting from the category of form rings to the category of abelian groups.

\begin{definition}
	In the quadratic case, the kernel of the group homomorphism
	$${\rm KQ_1}^{\lambda}(R[X],\Lambda[X])\rightarrow {\rm KQ_1}^{\lambda}(R,\Lambda)$$
	induced from  the  form ring homomorphism $(R[X],\Lambda[X])\rightarrow (R,\Lambda); X\mapsto 0$ is denoted by 
		${\rm NKQ_1}^{\lambda}(R,\Lambda)$. 
\end{definition}
\begin{remark}
	Since the class of $\Lambda$-quadratic groups is a subclass of the class of quadratic groups, 
		the local-global principle holds for $\Lambda$-quadratic groups. We use this fact throughout the next section.
\end{remark}

\section{Absence of torsion in ${\rm NKQ_1}^{\lambda}(R, \lambda)$} 

In \cite{STI}, J. Stienstra showed that 
${\nk}(R)$ is a ${\w}(R)$-module, where ${\w}(R)$ is the ring of big Witt vectors (cf.\ \cite{BL1} and \cite{WEL1}).
Consequently, in (\cite{WEL}, \S 3), C.\ Weibel 
observed that if $k$ is a unit in a commutative local ring $R$, then ${\sk}(R[X])$ has no 
$k$-torsion. Note that if $R$ is a commutative local ring, then 
${\sk}(R[X])$ coincides with ${\nk}(R)$; indeed, if $R$ is a local ring, then
${\SL}_n(R)={\E}_n(R)$ for all $n>0$. Therefore, we may replace $\alpha(X)$ by 
$\alpha(X)\alpha(0)^{-1}$ and assume that $[\alpha(0)]=[{\I}]$. In \cite{RB1}, the first author 
extended Weibel's result for arbitrary associative rings. In this section we prove the analog 
result for $\lambda$-unitary Bass nil-groups, {\it viz.} ${\rm NK_1GQ}^{\lambda}(R, \lambda)$, 
where $(R,\lambda)$ is the form ring as introduced by A.\ Bak in \cite{Bak1}. The main 
ingredient for our proof is an analog of Higman linearization (for a subclass of 
Bak's unitary group) due to V. Kopeiko; cf.\! \cite{V}. 
For the general linear groups, Higman linearization (cf.\! \cite{RB1}) allows us 
to show that ${\nk}(R)$ has a unipotent representative. The same result is not true in general 
for the unitary nil-groups. Kopeiko's results in \cite{V}, \cite{V1} give a complete description 
of the elements of ${\rm NK_1GQ}^{\lambda}(R, \lambda)$ which have (unitary) unipotent representatives.

\begin{definition}
	For an associative unital ring $R$ we consider the truncated polynomial ring 
	$$R_t=\frac{R[X]}{(X^{t+1})}.$$
\end{definition}

\begin{lemma}[{cf. \!\cite{RB1}, Lemma 4.1}] \label{la3} 
	Let $P(X)\in R[X]$ be any polynomial.
	Then the following identity holds in the ring $R_t$: 
	\begin{equation*}
		(1+X^r P(X))=(1+X^rP(0))(1+X^{r+1}Q(X)),
	\end{equation*}
	where $r>0$ and  $Q(X)\in R[X]$, with $\deg(Q(X))< t-r$.
\end{lemma}
\begin{proof} Let us write $P(X)=a_0+a_1X+\cdots+a_{t}X^{t}$. Then we can 
express the polynomial $P(X)$ as $P(X)=P(0)+XP'(X)$ for some $P'(X)\in R[X]$. Now, in $R_t$
{
	\begin{align*}
		(1+X^r P(X))(1+X^r P(0))^{-1} 
		& =  (1+X^r P(0)+X^{r+1}P'(X))(1+X^r P(0))^{-1}\\
		& = 1+X^{r+1}P'(X)(1-X^rP(0)+X^{2r}(P(0))^2-\cdots)\\
		& =  1+X^{r+1}Q(X),
\end{align*}}
\!\!where $Q(X)\in R[X]$ with $\deg(Q(X))< t-r$. 
Hence the lemma follows. 
\end{proof}

\begin{remark}
    Iterating the above process we can write, for any polynomial 
$P(X)\in R[X]$,
$$(1+XP(X))=\Pi_{i=1}^t(1+a_iX^i)$$ in $R_t$, 
for some $a_i\in R$. By ascending induction, it follows that the $a_i$'s 
are uniquely determined.
\end{remark}

\begin{lemma}
	\label{keylemma}
	Let $R$ be an associative ring, and $k\in \mbb{Z}$ such that $kR=R$. Let $P(X)$ be a polynomial in $R[X]$. Assume $P(0)$ lies in the center of $R$. Then, if
	$$(1+X^rP(X))^{k^r}=1 $$
 for some $r\geq 0$, we have that
 $$(1+X^rP(X))=(1+X^{r+1}Q(X))$$ in the ring $R_t$, for some $Q(X)\in R[X]$ with $\deg(Q(X))<t-r$.
\end{lemma}

The following result is due to V. Kopeiko, cf. \!\cite{V}. This is an analog of Higman linearization for this special case.

\begin{proposition} 
	\label{kopeiko NK_1} Let $(R, \lambda)$ be a form ring. Then, 
	every element of the group ${\rm NKQ_1}^{\lambda}(R,\Lambda)$ has a 
	representative of the form $$[a;b,c]_n=\begin{pmatrix}
		{\I}_{r}-aX & bX\\-cX^{n}& {\I}_r+a^*X+\cdots +(a^*)^nX^n
	\end{pmatrix}\in {\GQ}^{\lambda}(2r,R[X], \Lambda[X])$$ 
	for some positive integers $r$ and $n$, where $a,b,c\in {\rm M}_r(R)$ satisfy the following conditions:
	\begin{enumerate}
		\item the matrices $b$ and $ab$ are Hermitian and also $ab = ba^*$,
		\item the matrices $c$ and $ca$ are Hermitian and also $ca = a^*c$,
		\item  $bc = a^{n+1}$and $cb = (a^*)^{n+1}$.
	\end{enumerate}
\end{proposition}

\begin{corollary}
	\label{class is hyperbolic} Let $R$ be an associative ring.  
	Let $[\alpha]\in {\rm NKQ_1}^{\lambda}(R,\Lambda)$ have the representation  $[a;b,c]_n$ for some $a,b,c\in {\rm M}_n(R)$ according to Proposition  \ref{kopeiko NK_1}. 
	Then, if it is true that $({\rm I}_r-aX)\in {\GL}_r(R)$, we have $[\alpha]=[\mbb{H}({\rm I}_r-aX)]$ in ${\rm NKQ_1}^{\lambda}(R,\Lambda)$.
\end{corollary}

\begin{proof} By Lemma \ref{ifaisinvertible}, we have $[a;b,c]_n\equiv \mbb{H}({\I}_r-aX) \pmod{{\EQ}^{\lambda}(2r, R[X],\Lambda[X])}$. 
Hence  $[\alpha]=[\mbb{H}({\I}_r-aX)]$ in ${\rm NKQ_1}^{\lambda}(R, \Lambda)$. \end{proof}  

\begin{theorem} \label{th1} 
	\label{has no torsion} 
	Let $(R,\Lambda)$ be a form ring, where $R$ is an associative ring with $1$
	and $k$ is an invertible integer, i.e., $kR=R$. Let 
	$$[\alpha(X)]=\left[\begin{pmatrix}
		A(X)&B(X)\\C(X)&D(X)
	\end{pmatrix}\right]\in {\rm NKQ_1}^{\lambda}(R,\Lambda)$$ with 
	$A(X)\in {\GL}_r(R[X])$ for some $r\in \mathbb{N}$. 
	Then $[\alpha(X)]$ has no $k$-torsion.
\end{theorem}
\begin{proof}  By Theorem \ref{kopeiko NK_1}, $[\alpha(X)]=[[a;b,c]_n]$ for some $a,b,c\in {\rm M}_s(R)$ and for some natural numbers $n$ and $s$.  
Note that in Step $1$ of the Proposition \ref{kopeiko NK_1}, the invertibility of the top left corner of the matrix $\alpha$ is not 
changed during the linearization process.
The same is true for the remaining steps of the Proposition \ref{kopeiko NK_1}, and so, since the top left corner matrix is $A(X)\in {\GL}_r(R[X])$, we get $({\I}_s-aX)\in {\GL}_s(R[X])$. 
Using Corollary \ref{class is hyperbolic}, one gets $[\alpha(X)]=[\mbb{H}({\I}_s-aX)]$. Now let $[\alpha]$ have $k$-torsion, which implies that also $[\mbb{H}({\I}_r-aX)]$ has $k$-torsion. Since $({\I}_r-aX)$ is invertible, it follows that $a$ is nilpotent. 
Let $t$ be such that $a^{t+1}=0$. Since $[({\I}_r-aX)]^k=[{\I}]$ in ${\rm NKQ}_1^{\lambda}(R[X], \Lambda[X])$, 
by arguing as in \cite{RB}, we get $[{\I}_r-aX]=[I]$ in ${\rm NKQ}_1^{\lambda}(R[X], \Lambda[X])$. 
This completes the proof. \end{proof}  

	
	\section{Graded Analog}
	We recall the 
	well-known ``{\bf Swan--Weibel homotopy trick}'', which is the main ingredient to handle the graded case. Let 
	$R = R_0\oplus R_1 \oplus R_2 \oplus \cdots$ be a graded ring. 
	An element $a \in R$ is denoted by $a = a_0 + a_1 + a_2 + \cdots $, where $a_i\in R_i$ for each $i$, 
	and all but finitely many $a_i$'s are zero. Let $R_+= R_1 \oplus R_2 \oplus \cdots$. 
	A graded structure of $R$ induces a graded structure on ${\rm M}_n(R)$ (the ring of $n \times n$ matrices with entries in $R$).

	\begin{definition}
		Let $a \in R_0$ be a fixed element. We fix an element $b = b_0 + b_1 + \cdots$ in $R$ and define a ring homomorphism 
			$\epsilon: R \rightarrow R[X]$ as follows: \[ \epsilon(b) = \epsilon (b_0 + b_1 + \cdots )\; = \; b_0 + b_1X + b_2X^2 + \cdots + b_iX^i + \cdots .\]
			Then we evaluate the polynomial $\epsilon(b)(X)$ at $X = a$ and denote the image by $b^+(a)$ i.e., $b^+(a) = \epsilon(b)(a)$.
			Note that $\big(b^+(x)\big)^+(y) = b^+(xy)$. Observe, $b_0=b^{+}(0)$. 
			We shall use this fact frequently.
	\end{definition}

	The above ring homomorphism $\epsilon$ induces a group homomorphism from ${\GL}(n,R)$  to ${\GL}(n,R[X])$ of rank $n$ for every $n \geq 1$, i.e., for 
	$\alpha \in {\GL}(n,R)$ we get a map $$\epsilon: {\GL}(n,R) \rightarrow {\GL}(n,R[X]) \text{, defined by} $$ 
	$$\alpha = \alpha_0 \oplus \alpha_1 \oplus \alpha_2 \oplus \cdots \mapsto \alpha_0 \oplus \alpha_1X \oplus \alpha_2X^2 \cdots,$$
	where $\alpha_i\in {\rm M}(n,R_i)$. 
	As above, for $a \in R_0$, we define $\alpha^+(a)$ as $\alpha^+(a) = \epsilon(\alpha)(a).$

	The graded dilation lemma  and graded local-global principle are proved in \cite{BS} for linear, symplectic and orthogonal groups. 
	Arguing in similar manner one gets:
	
	%
	%
	\begin{theorem}\label{LG principle} {\bf (Graded Local-Global Principle)}
		Let $R = R_0\oplus R_1 \oplus R_2 \oplus \cdots$ be an almost commutative graded ring with identity $1$. 
		Let $\alpha\in {\GQ}(2n, R, \Lambda)$ be such that $\alpha\equiv {\I}_{2n}\pmod{R_+}$. If 
		$\alpha_{\mathfrak m}\in {\EQ}(2n, R_{\mathfrak m},\Lambda_{\mathfrak m})$, 
		for every maximal ideal $\mathfrak m\in \Max(C(R_0))$, then $\alpha\in {\EQ}(2n,R, \Lambda).$
	\end{theorem}
	
	Moreover, the L-G principle for the elementary subgroups and their normality property are equivalent.
	
	\begin{theorem} 
		
		Let $R = R_0\oplus R_1 \oplus R_2 \oplus \cdots$ be an almost commutative graded ring with identity $1$. 
		Then the following assertions are equivalent for $2n\ge 6$:
		\begin{enumerate}
			\item ${\EQ}(2n,R,\Lambda)$ is a normal subgroup of ${\GQ}(2n, R,\Lambda)$.
			\item If $\alpha\in {\GQ}(2n,R,\Lambda)$ with $\alpha^+(0)={\rm I}_{2n}$, and $\alpha_{\mathfrak{m}}\in {\EQ}(2n, R_{\mathfrak{m}}, 
			\Lambda_{\mathfrak{m}})$ for every maximal ideal $\mathfrak{m}\in {\rm Max}(C(R_0))$, then $\alpha\in {\EQ}(2n, R, \Lambda)$.
		\end{enumerate}
		
	\end{theorem}

	As an application of Theorem \ref{LG principle} and the Theorem \ref{th1}, we obtain the following:
	
	\begin{theorem} \label{th3} 
		Let $R =R_0\oplus R_1\oplus \dots$ be an almost commutative graded ring with $1$. 
		Let $N = N_0 +N_1 +\dots+N_r\in {\rm M}_r(R)$ be a nilpotent matrix, and let ${\I}$ denote the identity matrix. 
		Let $k\in \mbb{Z}$ be a unit in $R_0$. If
		$[({\I}+N)]^k=[{\I}]$ in ${\rm NKQ}_1^{\lambda}(R, \Lambda)$, then $[{\I} +N] = [{\I} +N_0]$.
	\end{theorem}
	\begin{proof} Consider the ring homomorphism $f:R\rightarrow R[X]$ defined by $$f(a_0+a_1+\dots)=a_0+a_1X+\dots.$$ Then
	\begin{align*}
		[({\I}+N)^k]=[{\I}]&\Rightarrow f([{\I}+N]^k)=[f({\I}+N)]^k=[{\I}]\\
		&\Rightarrow [({\I}+N_0+N_1X+\dots +N_rX^r)]^k=[{\I}].
	\end{align*}
	Let $\mathfrak{m}$ be a maximal ideal in $C(R_0)$.  By Theorem \ref{has no torsion}, 
	we have $$[({\I}+N_0+N_1X+\dots+N_rX^r)]=[{\I}]$$ in ${\rm NKQ}_1^{\lambda}(R_{\mathfrak{m}},\Lambda_{\mathfrak m})$. 
	Hence by using the local-global principle we conclude $$[({\I}+N)]=[{\I}+N_0]$$ 
	in ${\nk}{\GQ}^{\lambda}(R,\Lambda)$, as required. \end{proof}
	
	\subsection*{Acknowledgment}
	We thank the unknown referee for useful suggestions to improve  the manuscript. We also thank V. Kopeiko and Sergey Sinchuk for useful discussions on Bass's nil group. 
 Research by the first author was supported by SERB-MATRICS grant for the financial
year 2020--2021. Research by the second author was supported by IISER (Pune) post-doctoral research grant, 2020--2021.


\EditInfo{July 28, 2022}{February 9, 2023}{Roozbeh Hazrat}

\end{paper}